\documentclass[a4paper,11pt]{article}
\usepackage{xypic}
\usepackage{tikz}
\usepackage{amsmath,amsthm}
\usepackage{amssymb}
\usepackage{mathrsfs}
\def\diam{\mathop{\rm diam}}

\def\lim{\mathop{\rm lim}}

\def\phi{\varphi}

\def\pf{{\it Proof:}~}

\newtheorem{theorem}{Theorem}
\newtheorem{lemma}[theorem]{Lemma}
\newtheorem{proposition}[theorem]{Proposition}
\newtheorem{definition}[theorem]{Definition}
\newtheorem{corollary}[theorem]{Corollary}
\newtheorem{problem}[theorem]{Problem}
\newtheorem{remark}[theorem]{Remark}
\newtheorem{claim}[theorem]{Claim}
\newtheorem{examples}[theorem]{Examples}
\newtheorem{question}[theorem]{Question}
\newtheorem{conjecture}[theorem]{Conjecture}

\newcommand{\begintheorem}{\begin{theorem}}
\newcommand{\beginlemma}{\begin{lemma}}
\newcommand{\beginexercise}{\begin{exercise}}
\newcommand{\beginexercises}{\begin{exercises}}
\newcommand{\beginproposition}{\begin{proposition}}
\newcommand{\begindefinition}{\begin{definition}}
\newcommand{\begincorollary}{\begin{corollary}}
\newcommand{\beginproblem}{\begin{problem}}
\newcommand{\beginremark}{\begin{remark}}
\newcommand{\beginclaim}{\begin{claim}}
\newcommand{\beginassumptions}{\begin{assumptions}}
\newcommand{\beginexamples}{\begin{examples}}
\newcommand{\beginquestion}{\begin{question}}
\newcommand{\beginsassumptions}{\begin{sassumptions}}
\newcommand{\beginsassumption}{\begin{sassumption}}
\newcommand{\beginconjecture}{\begin{conjecture}}
\newcommand{\ud}{\, \mathrm{d}}

\providecommand{\abs}[1]{\lvert#1\rvert}
\providecommand{\norm}[1]{\lVert#1\rVert}
\providecommand{\Bnorm}[1]{\Bigl\lVert#1\Bigr\rVert}

\title{\bf Quenched invariance principle for random walks in balanced random environment}
\author{Xiaoqin Guo\thanks{School of Mathematics, University of Minnesota,
206 Church St SE, Minneapolis, MN 55455. Partially supported by NSF
grant  DMS-0804133.}
\and Ofer Zeitouni\thanks{School of Mathematics, University of Minnesota,
206 Church St SE, Minneapolis, MN 55455 and Faculty of Mathematics,
Weizmann Institute, Rehovot 76100, Israel. Partially supported by NSF
grant  DMS-0804133,  the Israel Science Foundation and the 
Herman P. Taubman chair of Mathematics at the Weizmann Institute}
}
\date{April 27, 2010. Revised August 23, 2011}
\begin{document}
\maketitle

\begin{abstract}
We consider random walks in a balanced random environment in
$\mathbb{Z}^d$, $d\geq 2$.
We first
prove an invariance principle (for $d\ge2$) 
and the transience of the random walks
when $d\ge 3$ (recurrence when  $d=2$)
in an ergodic environment 
which is not uniformly elliptic but satisfies certain moment condition. 
Then, using percolation arguments, we show that under mere ellipticity, 
the above results 
hold for random walks in i.i.d. balanced environments.
\end{abstract}

\section{Introduction}\label{introduction}
In recent years, there has been much interest in the study
of invariance principles and transience/recurrence
for random walks in random environments (on the
$d$-dimensional lattice $\mathbb{Z}^d$)
with non uniformly
elliptic transitions probabilities. Much of this work has 
been in the context of reversible models, either for walks on percolation
clusters or for the random conductance model, see
\cite{Bar04,SS05,MR05,BB,MaP07,Ma08,BarDe10}.  In those cases,
the main issue is the transfer of annealed estimates (given e.g.
in \cite{DFGW89}) to the quenched setting, and the control
of the quenched mean displacement of the walk. 
On the other hand, in these models the reversibility of the walk
provides for explicit expressions for certain invariant measures
for the environment viewed from the point of view of the particle.

The non-reversible setup has proved to provide many additional, and
at this point insurmountable, 
challenges, even in the uniformly elliptic setup, see 
\cite{Zrev} for a recent account, and it is therefore premature
to study in that generality the effects of non uniformly elliptic transition
probabilities. 
However, a particular class
for which the (quenched) invariance principle has 
been established in the uniformly elliptic 
setup is that of walks in balanced environments, see \cite{La}.
In that case, a-priori estimates of the Alexandrov-Bakelman-Pucci
type give enough control that allows one to prove the existence
of invariant measures (for the environment viewed from
the point of view of the particle), and the fact that the walk
is a (quenched) martingale together with ergodic arguments yield
the invariance principle (obviously, control of the quenched 
mean displacement,
which vanishes, is automatic). The establishment of recurrence (for $d=2$)
and transience (for $d\geq 3$) requires some additional
arguments, due to Kesten and Lawler, respectively, see
\cite{ZO} for details. 

It is our goal in this paper to explore the extent to which the
assumption of uniform ellipticity can be dropped in
this non-reversible, but balanced, setup. Not surprisingly,
it turns out that
some moment assumptions on the
 ellipticity constant suffice to yield the 
invariance principle in the ergodic environment setup, 
after some analytical effort has been expanded in obtaining a-priori 
estimates. What is  maybe more surprising is that for i.i.d.
environments, no assumptions of uniform ellipticity are needed at
all.

We describe now precisely the model 
we consider.
%consider a random walk in a balanced random environment on the 
%$d$-dimensional integer lattice $\mathbb{Z}^d$, $d\ge 2$. 
%The model is described as follows. 
Let $\mathcal{M}$ be the space of all probability measures on 
$V=\{v\in\mathbb{Z}^d: |v|\le 1\}$,
 %GGG02/15/2010 added "where $|\cdot|$ denotes the $l^2$-norm"
 where $|\cdot|$ denotes the $l^2$-norm. We equip  $\mathcal{M}$ with the weak topology on probability measures, which makes it into a Polish space, and equip $\Omega=\mathcal{M}^{\mathbb{Z}^d}$ with the induced Polish structure. Let $\mathcal{F}$ be the Borel
$\sigma$-field of $\Omega$ and $P$ a probability measure on $\mathcal{F}$.

A random \textit{environment} is  an element $\omega
=\{\omega(x, v)\}_{x\in{\mathbb{Z}^d}, v\in V}$ of $\Omega$ with distribution $P$. The random environment is called i.i.d. if $\{\omega(x, \cdot)\}_{x\in\mathbb{Z}^d}$ are i.i.d. across the sites $x$ under $P$. The random environment is called \emph{balanced} if $$P\{\omega(x, e_i)=\omega(x,-e_i) \mbox{ for all $i$ and all $x\in\mathbb{Z}^d$}\}=1,$$
and \textit{elliptic} if $P\{\omega(x,e)>0 \mbox{ for all $|e|=1$ and all $x\in\mathbb{Z}^d$}\}=1$.

The random walk in the
random environment $\omega\in\Omega$ (RWRE)
started at $x$ is the canonical Markov
chain $\{X_n\}$ on $(\mathbb{Z}^d)^\mathbb{N}$,
with state space $\mathbb{Z}^d$ and law $P_\omega^x$ specified by
\begin{align*}
&P_\omega^x\{X_0=x\}=1,\\
&P_\omega^x\{X_{n+1}=y+v | X_n=y\}=\omega(y, v), \quad v\in V.
\end{align*}
The probability distribution $P_\omega^x$ on $((\mathbb{Z}^d)^\mathbb{N}, \mathcal{G})$ is called the \textit{quenched law}, where $\mathcal{G}$ is the $\sigma$-field generated by cylinder functions.
Note that for each $G\in\mathcal{G}$, $P_\omega^x (G) : \Omega\to [0,1]$ is a $\mathcal{F}$-measurable function.
The joint probability distribution $\mathbb{P}^x$ on $\mathcal{F}\times\mathcal{G}$:
\[
\mathbb{P}^x (F\times G)=\int_{F} P_\omega^x (G)P(\ud\omega), \qquad F\in\mathcal{F},\, G\in\mathcal{G},
\]
is called the \textit{annealed} (or \textit{averaged})
law. Expectations with respect to $P_\omega^x$ and $\mathbb{P}^x$ are denoted
by $E_\omega^x$ and $\mathbb{E}^x$, respectively.

Define the canonical shifts $\{\theta^y\}_{y\in\mathbb{Z}^d}$ on
$(\Omega, \mathcal{F})$ by $(\theta^y \omega)(x,v)=\omega(x+y,v)$.
Throughout
the paper,
we always assume that the system $(\Omega, \mathcal{F}, P)$ is
ergodic with respect to the group of shifts $\{\theta^y\}$ and
that the environment is balanced and elliptic.\\

Let $o=(0,\cdots, 0)$ denote the origin and
\[
X_t^n :=\frac{1}{\sqrt{n}}X_{\lfloor tn\rfloor}+
                 \frac{tn-\lfloor tn\rfloor}{\sqrt{n}}(X_{\lfloor tn\rfloor+1}-X_{\lfloor tn\rfloor})
                                     , \quad t\ge 0.
\] We say that the \textit{quenched invariance principle} holds with
nondegenerate covariances if for $P$-almost every $\omega\in\Omega$,
the $P_\omega^o$ law of the path $\{ X_t^n\}_{t\geq 0}$ converges weakly
 to a Brownian motion on $\mathbb{R}^d$ with
covariance matrix $(a_i\delta_{ij})_{1\le i,j\le d}$, $a_i>0$,
as $n\to\infty$.

Lawler \cite{La} proved a quenched invariance principle
%{\bf Check!}
for random walks in a
balanced random environment
%is proved by Lawler
%$\cite{La}$
%(see detailed proofs in $\cite{Sz,ZO}$),
under the assumption that the random environment is {\it
uniformly elliptic}, i.e.
$$P\{\omega(x, e)\ge \varepsilon_0 \mbox{
for all $|e|=1$}\}=1 \mbox{\qquad for some $\varepsilon_0>0$.}$$

As mentioned above,
our goal in this paper is to study the extent to which the
uniform ellipticity assumption can be dropped.
Let
\begin{equation}
	\label{epsdef}
	\varepsilon(x)=\varepsilon_{\omega}(x):=
	[\prod_{i=1}^{d}\omega(x,e_i)]^{\frac{1}{d}}.\end{equation}
Our first main result
%, see Theorem  \ref{CLT1},
is that
if $\mathrm{E}\varepsilon(o)^{-p}< \infty$ for some $p>d$, then
the quenched invariance principle
holds and moreover, the RWRE is transient $P$-almost surely if
$d\geq 3$. (Recurrence for $d=2$ under 
the condition $E\varepsilon(0)^{-p}<\infty$
follows from the quenched invariance principle and ergodicity by an 
%OOO0828
unpublished argument of Kesten 
%GGG0827
detailed in
\cite[Page 281]{ZO}. Note that this argument cannot be used to
prove transience in dimensions $d\geq 3$, 
even given an invariance principle, since in higher dimensions
the invariance principle does not give useful
information on the range of the
random walk; the behavior of the
range is a crucial element in Kesten's argument.)

\begin{theorem}\label{CLT1}
%OOO1
Assume that the random environment is ergodic, elliptic
 and balanced.
\begin{enumerate}
\item[(i)] If $E\varepsilon(o)^{-p}< \infty$ for some $p>d\ge 2$,
%OOO
then
 the quenched invariance principle holds with a
nondegenerate limiting covariance.
\item[(ii)] If $E[(1-\omega(o,o))/\varepsilon(o)]^q< \infty$ for some $q>2$ and $d\ge 3$, then the RWRE is transient $P$-almost surely.
%GGG 3/11/2010
\end{enumerate}
\end{theorem}
\noindent
%OOO1
That some integrability condition on the tail of $\varepsilon(o)$ is needed
for part (i) to hold
is made clear by the (non-Gaussian) scaling limits of random walks in
Bouchaud's trap model, see \cite{Bou,BAC}.
%see e.g. \cite{BAC} and references therein.
In fact,
it follows from that example that Theorem \ref{CLT1}(i), or even an annealed
version of the CLT,
cannot
%OOO1
hold in general with  $p<1$.
%is easy to construct an elliptic example with

 The proof of Theorem \ref{CLT1}
 is based on a sharpening of
the arguments in
\cite{La,Sz,ZO}; in particular, refined versions of the maximum
principle for walks in balanced environments (Theorem \ref{MP})
and of a mean value inequality (Theorem
\ref{mvi})
play a crucial role.

When the environment is i.i.d. and elliptic, our second main result
%(Theorem \ref{CLT2})
is that
%OOO when the environment is elliptic and
if $|X_{n+1}-X_n|=1$ a.s., then the quenched invariance principle holds.
%GGG 04/05/2010
Moreover, the RWRE is $P$-almost surely transient when $d\ge 3$.
The proofs combine percolation arguments with Theorem \ref{CLT1}.
\begin{theorem}\label{CLT2}
Assume that the random environment is i.i.d., elliptic and
balanced.
\begin{enumerate}
\item[(i)] If $P\{\max_{|e|=1}\omega(o,e)\ge \xi_0\}$=1 for some positive constant $\xi_0$, then the 
quenched invariance principle holds with a non-degenerate limiting covariance.
\item[(ii)] When $d\ge 3$, the RWRE is transient $P$-almost surely.
\end{enumerate}
\end{theorem}
Because the 
transience or recurrence of the random walks does not change
if one considers the walk restricted
to its jump times,
one concludes, using Kesten's argument and the invariance
principle, compare with
Theorem \ref{CLT1}, that
for $d=2$, a random walk in a balanced elliptic i.i.d. 
random environment is recurrent $P$-a.s.

%GGG
Our proof of the invariance principles, like that of \cite{La}, is based
on the approach of the
``environment viewed from the point of view of the particle".
Specifically, set $\bar{\omega}(n)=\theta^{\mathrm{X}_n}\omega$,
then the process $ \bar{\omega}(n) $ is a Markov chain under
$ \mathbb{P}^{o} $ with state space $ \Omega $ and transition kernel \[ M(\omega',\ud\omega)=\sum_{i=1}^{d}[\omega'(o,e_i)\delta_{\theta^{e_i}\omega'}+
\omega'(o,-e_i)\delta_{\theta^{-e_i}\omega'}]+\omega'(o,o)\delta_{\omega'}. \]

%{\bf Adapt text below to cover invariance principle}

%OOO
Since $\{X_n\}$ is a (quenched) martingale,
standard arguments (see the proof of Theorem 6.2 in \cite{BB}) show that
the quenched invariance principle holds
whenever  an invariant
measure $Q\sim P$ of $\{\bar{\omega}(n)\}$ exists.
 The
%OOO
approach
of Lawler \cite{La}, which is a discrete version of the argument of Papanicolaou and Varadhan \cite{PV}, is to construct such a measure as the limit of invariant measures of periodized environments. We will
 follow this strategy using, as in \cite{Sz,ZO}, variants of
 \cite{KT} to derive estimates on solutions of linear elliptic difference
 equations. In the i.i.d. setup of
%OOO1
Theorem \ref{CLT2}, percolation estimates are used to control
pockets of the environment where those estimates are not strong enough.\\
% . Our proof of part
% (i) of Theorem \ref{CLT1} is a
% modification of the arguments in
% $\cite{La,Sz,ZO}$.
 %, however, in the proof of Theorem \ref{CLT2} we need some tools in percolation.\\

%GGG 04/05/2010
For the proof of the transience in the ergodic case, 
we use a mean value inequality and follow \cite{ZO}.
To prove the transience in the iid case, we employ percolation
arguments together with
a new maximum principle (Theorem \ref{mp2}) for walks with (possibly)
big jumps.
% and percolation arguments.\\

The structure of this paper is as follows.
In Section \ref{SePeEn} we construct the ``periodized environments" as in \cite{Sz, ZO}, and show that the
proof of $Q\sim P$ can be reduced to the proof of
%OOO
the inequality (\ref{Phi0}).
Using the  maximum principle, we then prove
(\ref{Phi0})
in Section \ref{SeMP}
under the assumptions of Theorem \ref{CLT1}(i).
In
Section \ref{SePercE}, devoted to the i.i.d. setup,
we prove
Theorem \ref{CLT2}(i),
%(\ref{Phi0}) in the i.i.d. setup,
using percolation tools. Section \ref{SeTran} is devoted
to the proof of the transience of the RWRE for $d\geq 3$, thus providing
a proof of Theorem \ref{CLT1}(ii). 
%GGG 04/05/2010
In Section \ref{SeTriid}, we will show
a modified maximum principle for balanced difference operators, and use it to prove
Theorem \ref{CLT2}(ii).
%We conclude
%in
%Section \ref{counterexample}
%with a counterexample to the invariance principle (and even the
%CLT) when $P(\omega(o,o)>0)>0$
%and $\varepsilon(o)^{-1}$ does not possess a $2+\delta$ moment.

%OOO
Throughout the paper, $C$ denotes a generic positive constant, that may depend on
dimension only, and whose value may change from line to line.

\section {The periodized environments}\label{SePeEn}
As in \cite{Sz, ZO}, the following periodic structure of the environment
is introduced.

Let $ \Delta_N (x_0)=\{x\in \mathbb{Z}^{d}: |x-x_0|_{\infty}\le N\} $ be the cube centered at $x_0$ of length $2N$. Let $\Delta_N=\Delta_N(o)$. For any $ x\in\mathbb{Z}^{d} $, set
$$ \hat{x}:=x+(2N+1)\mathbb{Z}^{d}\in \mathbb{Z}^{d}/(2N+1)\mathbb{Z}^{d}. $$

For any fixed $ \omega \in \Omega $, we define $\omega^{N}$ by setting
%to be such that
$ \omega^{N}(x)=\omega(x) $ for $ x \in \Delta_{N} $
and $\omega^N (y)=\omega^N (x)$ for $y \in \mathbb{Z}^{d} $ whenever $\hat{y}=\hat{x}$. Let $ \Omega^{N}=\{\omega^{N}: \omega\in\Omega\}$. Let $ \{X_{n,N}\} $ denote the random walk on $ \mathbb{Z}^{d} $ in the environment $ \omega^{N} $. Then $ \{\hat{X}_{n,N}\} $ is an irreducible finite-state Markov chain,
hence it possesses a unique invariant probability measure, which can
always be written in
the form 
%GGG0827 
%OOO0828
\[ \dfrac{1}{(2N+1)^d}\sum_{x\in\Delta_{N}}\Phi_{N}(x)\delta_{\hat{x}} . \]
Here $\Phi_N$ is some function on $\Delta_N$ and $(2N+1)^{-d}\Phi_{N}(\cdot)$ 
sums to $1$, so that $\Phi_N$ can be interpreted as a density with
%GGG0829
respect to the uniform measure on $\Delta_N$.
%is 
%defined as the density.
%GGG0827

Define 
\[ Q_{N}=Q_{N,\omega}=\dfrac{1}{(2N+1)^{d}}\sum_{x\in\Delta_{N}}\Phi_N (x)\delta_{\theta^{x}\omega^{N}} \]
as a probability measure on $\Omega^N$.
Then for any $x\in \Delta_N$, 
%GGG0829
\begin{align*}
\sum_{y\in\Delta_N} Q_N(\theta^y \omega^N)M(\theta^y \omega^N, \theta^x\omega^N)
&=\sum_{y\in\Delta_N} \frac{\Phi_N (y)}{(2N+1)^d}\omega^N(y,x)\\
&= \frac{\Phi_N (x)}{(2N+1)^d}=Q_N (\theta^x \omega^N).
\end{align*}
This implies that
$Q_N$ is the invariance probability measure
%OOO0827
(with respect to the kernel $M$) for the Markov chain 
$ \{ \bar{\omega}^{N}(n)\} $ on $\Omega^{N}  $.

We will show that $Q_N $ converges weakly to some measure $ Q $ with good properties. To do this, we first introduce a sequence of measures
\[ P_{N}=P_{N,\omega}=\dfrac{1}{(2N+1)^{d}}  \sum_{x\in\Delta_{N}}\delta_{\theta^{x}\omega^{N}} \]
which by the multidimensional ergodic theorem (see Theorem (14.A8) in \cite{Ge} and also Theorem 1.7.5 in \cite{Kr}) converges weakly to $ P $, $P\mbox{- }a.s.  $

Let $ \{\omega_{\gamma}^{N}\}_{\gamma=1}^{k} $ denote the set of distinct states in $ \{ \theta^{x}\omega^{N}\}_{x\in \Delta_{N}} $ and $ C_{N}(\gamma):=\{x\in \Delta_{N}: \theta^{x}\omega^{N}=\omega_{\gamma}^{N}\} $.
Set,
for any finite subset $E\subset\mathbb{Z}^d$,
$$ \lVert f\rVert_{E,j}:=(|E|^{-1}\sum_{x \in E} |f(x)|^{j})^{\frac{1}{j}}.$$
Since
$ \ud Q_N/\ud P_N=\sum_{\gamma=1}^{k}\delta_{\omega_{\gamma}^{N}}|C_{N}(\gamma)|^{-1}\sum_{x\in C_{N}(\gamma)}\Phi_{N}(x):=f_{N}$, we have that for any measurable function $g$ on $\Omega$,
\begin{align}\label{Q_n0}
|Q_N g|
 &\le (\int  f_{N}^{\alpha} \ud P_N)^{\frac{1}{\alpha}}(\int  |g|^{\alpha'} \ud P_N)^{\frac{1}{\alpha'}} \nonumber\\
 &\le  \big( \frac{1}{|\Delta_N |}\sum_{\gamma =1}^{k}\sum_{x \in C_N (\gamma)}\Phi_N (x)^{\alpha} \big)^{\frac{1}{\alpha}}(\int  |g|^{\alpha'} \ud P_N)^{\frac{1}{\alpha'}}
 \nonumber\\
 &=  \lVert         \Phi _N \rVert _{\Delta_N, \alpha}(P_{N}|g|^{\alpha'})^{\frac{1}{\alpha'}},
\end{align}
where
$\alpha'$ is the H\"older conjugate of $\alpha$, $1/\alpha+1/\alpha'=1$,
and we used H\"older's inequality in the first
and the second inequalities.
%and
%$\alpha'$ is the H\"older conjugate of $\alpha$, $1/\alpha+1/\alpha'=1$.
Since $\Omega$ is compact with respect to the product topology,
along some subsequence $N_k\to\infty$,
$\{Q_{N_k}\}$ converges weakly to a limit, denoted $Q$.
%of $ \{Q_N\} $.
Assume for the moment that
\begin{equation}\label{Phi0}
\varlimsup_{N\to\infty}\lVert         \Phi _N \rVert _{\Delta_N, \alpha}
\le C, \quad P\mbox{- }a.s..
\end{equation}
We show that then, for a.e. $\omega\in\Omega$,
\begin{equation}\label{f2}
Q\ll P   .
\end{equation}
Indeed, let $A\subset \Omega$ be measurable.
Let $ \rho $ denote a metric on the Polish space $ \Omega $.
For any closed subset
$ F \subset A $, $ \delta >0 $,
introduce the function $f(\omega)=[1-\rho (\omega, F)/\delta ]^+$
which is supported on
$ F_{\delta}=\{\omega \in \Omega: \rho (\omega, F) <\delta\}$.
%Here $ \rho $ is a metric of the Polish space $ \Omega $.
Then by (\ref{Q_n0}), (\ref{Phi0}),
\begin{equation*}
Q F\le \varlimsup_{N\to\infty}Q_N f \le C (P f^{\alpha'})^{\frac{1}{\alpha'}}\le  C (P F_{\delta})^{\frac{1}{\alpha'}} .
\end{equation*}
Letting  $ \delta\downarrow 0 $, we get $ Q F \le C (P F)^{\frac{1}{\alpha'}}$.
Taking supremums over all closed subset $ F \subset A $, one concludes that
$Q A \le C\cdot (P A)^{\frac{1}{\alpha'}}$, which proves (\ref{f2}).

Once we have (\ref{f2}), it is standard to check, using ellipticity,
that $ \bar{\omega}(n) $ is ergodic with respect to $ Q $ and $Q\sim P$
(see \cite{Sz, ZO}). (Thus,
by the ergodic theorem, $Q$ is uniquely determined by $Qg=\lim_{n\to\infty} E\sum_{j=0}^{n-1}g(\bar{\omega}_j)/n$ for every bounded measurable $g$. Hence $Q$ is \textit{the} weak limit of $Q_N$.) Therefore, to prove the invariance principle it suffices to prove (\ref{Phi0}).
Sections
\ref{SeMP} and Section \ref{SePercE} are devoted to the proof of (\ref{Phi0}),
under the assumptions of Theorems \ref{CLT1} and \ref{CLT2}.

%OOO
\section {Maximum Principle and proof of Theorem \ref{CLT1}(i)}\label{SeMP}
Throughout this section, we fix an $\omega\in \Omega$.
For any bounded set $ E \subset \mathbb{Z}^{d} $, let $\partial E =\{y \in E^{c}: \exists x\in E, |x-y| _{\infty}=1\}$, $ \bar{E}=E \bigcup \partial E $ and $ \diam(E)=\max\{|x-y| _{\infty}: x, y\in E\} $. For any function $f$ defined on $ \bar{E}$ ,
let $ L_{\omega} $ denote the operator
\[ (L_{\omega}f )(x)=\sum_{i=1}^{d}\omega(x, e_i)[f(x+e_i)+f(x-e_i)-2f(x)], \quad x\in E .\]

The following discrete maximum principle is an adaption of Theorem 2.1 of \cite{KT}.
%GGG
%The following theorem is a modification of
%Theorem 2.1 of $\cite{KT}$, see also $\cite{La, ZO}$.
\begin{theorem}[Maximum Principle]\label{MP}
Let $E\subset \mathbb{Z}^d$ be bounded,
and  let $u$ be a function on $ \bar{E} $. For all
$x\in E$, assume
$ \varepsilon(x)>0 $ and define
$$ I_{u}(x):=\{s \in \mathbb{R}^{d}:  u(x)-s\cdot x \ge u(z)-s\cdot z,  \forall z \in\bar{E}\}  .$$
If $ L_{\omega} u(x) \ge -g(x)$ for all $x \in E$ such that $I_u (x)\ne \emptyset$,
then
\begin{equation}\label{mpremark}
 \max_{E} u \le
%GGG \dfrac{1}{2}
  C\diam{\bar{E}}
%GGG \cdot |E|^{\frac{1}{d}}
 \bigg(\sum_{\substack{
                           x\in E\\
                           I_u (x)\ne \emptyset}
                      }|\frac{g}{\varepsilon}|^d\bigg)^{\frac{1}{d}}
 +\max_{\partial E}u  .
\end{equation}
In particular,  \[ \max_{E} u \le
%GGG \dfrac{1}{2}
C\diam{\bar{E}}\cdot |E|^{\frac{1}{d}} \lVert         \frac{g}{\varepsilon}\rVert _{E,d}+\max_{\partial E}u   .\]
\end{theorem}
%GGG 04/05/2010
\pf See the proof of Theorem 2.1 in \cite{KT}.\qed\\
%GGG: The proof of the maximum principle is removed.
%\pf
%Without loss of generality, assume $ g \ge 0 $ and \[ \max_{\bar{E}} u = u (x_0) > \max_{\partial E}u \]
%for some $ x_0 \in E $. Otherwise, there is nothing to prove.
%
%For $ s\in \mathbb{R}^{d} $ such that $ |s|_{\infty} \le [u(x_0)-\max_{\partial E}u]/(d\diam \bar{E}) $
%we have $ u(x_0)-u(x) \ge s \cdot (x_0-x) $ for all $ x \in \partial E $, which implies that $ s \in \bigcup_{x \in E} I_{u}(x) $.
%Hence
%
%\begin{equation}\label{f4}
% [-\dfrac{u(x_0)-\max_{\partial E}u}{d\diam \bar{E}}, \dfrac{u(x_0)-\max_{\partial E}u}{d\diam \bar{E}}]^{d} \subset \bigcup_{x\in E} I_{u}(x)  .
% \end{equation}
%
%Further, for $ s \in I_{u}(x) $ we have $u(x+e_i)-u(x)\le s\cdot e_i \le u(x)-u(x-e_i)  $,
%so
%
%\begin{equation}\label{f5}
% \bigcup_{x \in E} I_{u}(x) \subset \bigcup_{x \in E, I_{u}(x)\ne \emptyset} \prod_{i=1}^{d}[u(x+e_i)-u(x), u(x)-u(x-e_i)]
% \end{equation}
%
%By (\ref{f4}), (\ref{f5}) we have
%\begin{eqnarray*}
%&&\Bigl( \dfrac{2u(x_0)-2\max_{\partial E} u}{d\diam \bar{E}}\Bigr)^{d}\\
% &\le & \sum_{x\in E, I_{u}(x)\ne \emptyset}  \prod_{i=1}^{d}[2u(x)-u(x-e_i)-u(x+e_i)] \\
%&\le & \sum_{x\in E, I_{u}(x)\ne \emptyset}\Bigl ( \dfrac{1}{d}\sum_{i=1}^{d}\frac{\omega(x, e_i)}{\varepsilon(x)}[2u(x)-u(x-e_i)-u(x+e_i)] \Bigr )^{d}\\
%&\le &  \sum_{\substack{
%                           x\in E\\
%                           I_u (x)\ne \emptyset}
%                      }|\dfrac{g(x)}{d \varepsilon(x)}|^{d}.
%\end{eqnarray*}
%This completes the proof.
%\qed

Define the stopping times $ \tau_0=0 $, $ \tau_1 =\tau :=\min \{j \ge 1: |X_{j,N}-X_{0,N}|_{\infty}> N\} $ and
$ \tau_{j+1}=\min \{n>\tau_j : |X_{n,N} -X_{\tau_j, N}|_{\infty}> N\} $.

\begin{lemma}\label{E}
Let $\omega^{N}$, $\{X_{n,N}\}$ be as in Section 1 and $\tau$ as defined above, then
there exists a constant $c$ such that, for all $N$ large,
%OOO
 $$E_{\theta^{x}\omega^{N}}^{o} (1-\frac{c}{N^2})^{\tau}\leq C <1   .$$
\end{lemma}
\pf
Since P is balanced, $ X_{n,N} $ is a martingale
and it follows from Doob's inequality that for any $K \ge 1  $,
\begin{align*}
P_{\theta^{x}\omega^{N}}^{o} \{\tau \le K\}
& \le 
 2 \sum_{i=1}^{d} P_{\theta^{x}\omega^{N}}^{o}\{\sup_{n \le K} X_{n, N}(i) \ge N+1\}\\
& \le 
 \frac{2}{N+1} \sum_{i=1}^{d}E_{\theta^{x}\omega^{N}}^{o} X_{K, N} (i)^{+}
\le \frac{2d}{N+1} \sqrt{K},
\end{align*}
where $ X_{n,N} (i) $ is the $i$-th coordinate of $ X_{n,N} $. Hence
\begin{equation*}
E_{\theta^{x}\omega^{N}}^{o} (1-\frac{c}{N^2})^{\tau}
 \le  (1-\frac{c}{N^2})^{K} +\frac{2d}{N+1}\sqrt{K}   .
 \end{equation*}
Taking $c=16 d^2$ and $K=N^2/16d^2$, we get
$E_{\theta^{x}\omega^{N}}^{o} (1-\frac{c}{N^2})^{\tau}\le e^{-1}+2^{-1} .$ \qed

\begin{theorem}\label{Phi}
\begin{equation}\label{Phi1}
\lVert         \Phi_N \varepsilon\rVert _{\Delta_N , \beta}\le C,
\end{equation}
where $\beta=d'=d/(d-1)$.
\end{theorem}
\pf
 Let $c$ be the same
constant as in the previous lemma. For any function $ h\ge 0 $ on $\Delta_N $,
\begin{align*}
&\lVert         \Phi_N \cdot h\rVert _{\Delta_N , 1}\\
&=
 \frac{c}{N^2}\sum_{x \in \Delta_N}\frac{\Phi_N (x)}{|\Delta_N |}
\sum_{m\ge 0}E_{\omega^N}^{x}\sum_{\tau_m\le\ j <\tau_{m+1}}(1-\frac{c}{N^2})^{j}h(\hat{X}_{j,N})\\
&\le 
 \frac{c}{N^2}\sum_{x \in \Delta_N}\frac{\Phi_N (x)}{|\Delta_N |}\sum_{m\ge 0}E_{\omega^N}^{x}(1-\frac{c}{N^2})^{\tau_m}
E_{\omega^N}^{\hat{X}_{\tau_m , N}}\sum_{j=0}^{\tau -1}h(\hat{X}_{j,N})\\
&\le 
 \frac{c}{N^2}\sum_{x \in \Delta_N}\frac{\Phi_N (x)}{|\Delta_N |}\sum_{m\ge 0}\big[\sup_{y\in \Delta_N}E_{\omega^N}^{y}(1-\frac{c}{N^2})^{\tau}\big]^{m}\cdot \sup_{y\in \Delta_N}E_{\omega^N}^{y}
\sum_{j=0}^{\tau -1}h(\hat{X}_{j,N}).
 \end{align*}
%GGG0827
Since the function $f(x)=E_{\omega^N}^{x} \sum_{j=0}^{\tau -1}h(\hat{X}_{j,N})$ satisfies
\begin{equation}
\left\{
\begin{array}{rl}
L_{\omega^N}f(x)=h(x), & \text{if } x\in\Delta_N\\
f(x)=0, & \text{if } x\in\partial \Delta_N,
\end{array}
\right.
\end{equation}
we can apply the maximum principle (Theorem \ref{MP}) and get
$$\sup_{y\in \Delta_N}E_{\omega^N}^{y}\sum_{j=0}^{\tau -1}h(\hat{X}_{j,N}) \le C N^2
\lVert
\frac{h}{\varepsilon}
\rVert _{\Delta_N ,d}.$$
This together with Lemma \ref{E} and $\sum_{x \in \Delta_N}\Phi_N (x)/|\Delta_N |=1$ yield
$$\lVert         \Phi_N \cdot h\rVert _{\Delta_N , 1}
\le C \lVert         \frac{h}{\varepsilon}\rVert _{\Delta_N ,d}.$$
Hence by the duality of norms,
\begin{equation*}
\lVert         \Phi_N \varepsilon\rVert _{\Delta_N , \beta}=\sup_{\lVert         \frac{h}{\varepsilon}\rVert _{\Delta_N, d}=1}\lVert
\Phi_N h\rVert _{\Delta_N, 1}\le C .  \quad  \mbox{\qed}
\end{equation*}

\noindent{\it Proof of (\ref{Phi0}) under the assumption of Theorem \ref{CLT1}(i) :}\\
Assume that
\begin{equation}\label{asm}
\mathrm{E}\varepsilon(o)^{-p}< \infty
\mbox{ for some } p>d.
\end{equation}
Take $\alpha=(1-1/d+1/p)^{-1}$. We use H\"older's inequality
%OOO
and Theorem \ref{Phi}
to get
\[
\lVert \Phi_N\rVert _{\Delta_N,\alpha}\le \lVert  \Phi_N \varepsilon\rVert _{\Delta_N , \beta}\lVert         \varepsilon^{-1}\rVert _{\Delta_N, p}
\le C \lVert \varepsilon^{-1}\rVert _{\Delta_N, p}  .
\]
By the multidimensional ergodic theorem,
\begin{equation*}
%GGG04/06/2010	\label{eq-160110}
\lim_{N \to \infty}\lVert         \varepsilon^{-1}\rVert _{\Delta_N, p}=(E \varepsilon(o)^{-p})^{\frac{1}{p}}<\infty,
\quad  P\mbox{- }a.s. . \mbox{\qed }
\end{equation*}

\begin{remark}
Without the assumption (\ref{asm}),  the conclusion
%GGG \eqref{eq-160110}
\eqref{Phi0} may fail.
%things are more complicated.
To see the difficulty,
%when (\ref{asm}) does not hold,
%GGG
let $$A=A(\omega, \varepsilon_0)=\{x:\min_i \omega(x, e_i)<\varepsilon_0\}.$$
By (\ref{Phi1}) we have
$$\lVert         \Phi_N 1_{A^c}\rVert _{\Delta_N, \beta}\le \lVert         \Phi_N \frac{\varepsilon}{\varepsilon_0}\rVert _{\Delta_N, \beta} \le \frac{C}{\varepsilon_0}.$$
In order to proceed as before, we need to show that
$\varlimsup_{N\to\infty}\lVert \Phi_N 1_{A}\rVert _{\Delta_N, \alpha} \le C$
for some $1 <\alpha \le \beta$ .
As Bouchaud's trap model
\cite{Bou,BAC} shows,
%the counterexample in Section \ref{counterexample} shows,
this is not always the case. However, if $P\{\max_{|e|=1}\omega(o,e)\ge\xi_0\}=1$,
then for $x\in A$, we have, using that the environment is balanced,
some control of $\Phi_N (x)$
by $\Phi_N|_{A^c}$ (see Lemma \ref{Phicontrol}).
%which is an interesting feature of the balanced environment.
Further,
%Also note that
in the i.i.d. case, $A$ corresponds to  a `site percolation' model,
whose cluster sizes can be estimated.
We will show in the next section that these properties lead
to a proof of (\ref{Phi0}) in the i.i.d. setup, without moment assumptions.
\end{remark}

\section{A percolation estimate and proof of Theorem \ref{CLT2}(i)}\label{SePercE}
In this section we consider the RWRE in the i.i.d. setting where
%GGG04/13/01
$\max_{|e|=1}\omega(x,e)\ge \xi_0$  for all $x\in\mathbb{Z}^d$ and all $\omega\in\Omega$.
%By choosing $\varepsilon_0$ sufficiently small, we will prove (\ref{Phi0}) using a percolation estimate.\\
We begin by introducing some terminology.\\

The \textit{$l^1$-distance} (graph distance) from $x$ to $y$ is defined as
$$d(x,y)=|x-y|_1=\sum_{i=1}^{d}|x_i-y_i|.$$
Note that $|x|_{\infty}\le |x|_1 \le d|x|_{\infty}$.

In an environment $\omega$, we say that a site $x$  is \textit{open}(\textit{closed}) if $\min_i \omega(x, e_i)<\varepsilon_0 (\ge\varepsilon_0, resp.)$ and that an edge of $\mathbb{Z}^d$ is open if its endpoints are open. 
%GGG3/30/2010
Here $\varepsilon_0>0$ is a constant whose value is to be determined.
An edge is called closed if it is not open.
Let $A=A(\omega)$ denote the subgraph of $\mathbb{Z}^d$ obtained
by deleting all  closed edges and
closed sites.
We call $A(\omega)$ a \textit{site percolation}
with parameter $p=p(\varepsilon_0)=P\{\min_i \omega(x, e_i)< \varepsilon_0\}$.

A \textit{percolation cluster} is a connected component of  $A$.
(Although here a percolation cluster is defined as a graph, we also use it as
a
synonym for its set of vertices.)
The $l^1$ diameter of a percolation cluster $B$ is defined as $l(B)=\sup_{x\in B, y\in \partial B}d(x,y)$.
For $x\in A$, let $A_x$ denote the percolation cluster that contains $x$ and let $l_x$ denote its diameter. Set $A_x=\emptyset$ and $l_x=0$ if $x\notin A$.
%GGG3/30/2010
We let $\varepsilon_0$ be small enough such that $l_x<\infty$ for all $x\in \mathbb{Z}^d$.\\

%GGG3/30/2010
We call a sequence of sites $(x^1, \cdots, x^n)$ a \textit{path}  from $x$ to $y$ if $x^1=x$, $x^n=y$ and
$|x^j-x^{j+1}|=1$
%GGG08/25/2011 $d(x^j-x^{j+1})=1$ should be $|x^j-x^{j+1}|=1$
 for $j=1,\cdots, n-1$. Let
$$\square=\{(\kappa_1,\cdots, \kappa_d)\in \mathbb{Z}^d: \kappa_i=\pm 1\}.$$ 
We say that a path $\{x^1, \cdots, x^n\}$ is a \textit{$\kappa$-path}, $\kappa\in \square$,  if
$$\omega(x^j, x^{j+1}-x^j)\ge\xi_0$$ and
$\kappa_i(x^{j+1}-x^j)_i\ge 0$ for all $i=1,\cdots, d$ and $j=1,\cdots, n-1$.
%GGG3/30/2010
Observing that for each site there exist at least two neighbors (in opposite directions) to whom the transition probabilities are $\ge \xi_0$, we have the following property concerning the structure of the balanced environment:
\begin{itemize}
%GGG0827
\item  For any $x\in A$ and any $\kappa\in \square$, there exists a $\kappa$-path from $x$ to some $y\in\partial A_x$, and this path is contained in $\bar{A}_x$.
\end{itemize}

%GGG3/30/2010 deleted: We begin with a useful inequality.
This property gives us a useful inequality.
\begin{lemma}\label{Phicontrol}
For $x\in A\cap \Delta_N$, if $l_x \le N$, then
\begin{equation}\label{f7}
\Phi_N (x) \le \xi_0 ^{-l_x} \sum_{y \in \partial A_x \cap \Delta_N} \Phi_N (y).
\end{equation}
\end{lemma}
\pf Suppose that $A_x\neq\emptyset$ (otherwise the proof is trivial). Since $l_x \le N$,  $\bar{A}_x\subset \Delta_N (x)$. Note that at least one of the $2^d$ 
%GGG3/30/2010 vertices-->corners
corners of $\Delta_N (x)$ is contained in $\Delta_N$. Without loss of generality, suppose that
$v=x+(N,\cdots,N)\subset \Delta_N$. Then there is a $(1,\cdots, 1)$-path in $\bar{A}_x$ from $x$ to
some $y\in \partial A_x\cap \Delta_N$,
%The key observation is that starting from $x$, there is a path in $\prod_{i=1}^d [x_i, x_i+N]\cap \mathbb{Z}^d$ to some $y\in\partial A_x \cap \Delta_N$ such that the transition probability of each step is $\ge \xi_0$.  The choice of the path is explained as follows. Observe that  for each point there exists two neighbors (in opposite directions) to whom the transition probabilities are $\ge \xi_0$. We then move to the neighbor such that the $l^1$-distance from $x$ increases by $1$ and the $l^1$-distance to $v$ decreases by 1.
%This path will finally hit some $y\in\partial A_x \cap\Delta_N$, 
as illustrated in the following figure:

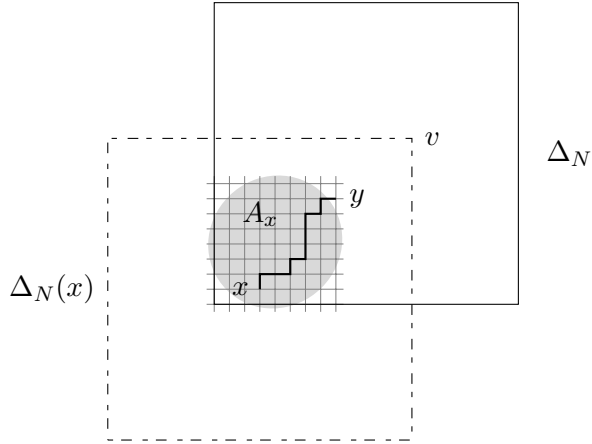
\begin{center}
\begin{tikzpicture}[scale=.2]{centered}
\fill[gray!30!white][rotate=45] (10, 1.5) ellipse (4.5 and 4.3);
\draw[step=1, gray, very thin] (1.5,3.5) grid (10.5, 12.5);
\draw[dash pattern=on 2pt off 3pt on 4pt off 4pt] (-5,-5) rectangle (15, 15) node[right=1pt]{$v$};
\draw (2,4) rectangle (22, 24);
\draw[thick] (5,5)--(5,6)--(7,6)--(7,7)--(8,7)--(8,10)--(9,10)--(9,11)--(10,11)node[right=1pt]{$y$};
\draw[left=1pt] (5,5) node{$x$};
\draw (3,10) node[right=1pt]{$A_x$};
\draw (23, 14) node[right=1pt]{$\Delta_N$};
\draw (-5, 5) node[left=1pt]{$\Delta_N(x)$};
\end{tikzpicture}
\end{center}

Recalling that $\Phi_N$ is the invariant measure for $\{\hat{X}_{n,N}\}$ defined in Section 1, we have
\begin{align*}
\Phi_N (y) 
&=
 \sum_{z\in \Delta_N} \Phi_N (z) P_{\omega^N}^{d(x,y)} (\hat{z}, \hat{y})\\
&\ge 
 \Phi_N (x) P_{\omega^N}^{d(x,y)}(\hat{x}, \hat{y})
\ge  \Phi_N (x) \xi_0^{l_x}.
\end{align*}
Here $P^m_{\omega^N} (\hat{z}, \hat{y})$ denotes the $m$-step transition probability
of $\{\hat{X}_{n,N}\}$ from $\hat{z}$ to $\hat{y}$.
\qed

Let $S_n=\{x: |x|_{\infty}=n\}$ denote the boundary of $\Delta_n$.
Let $x\to y$ be the event that $y\in \bar{A}_x$
and $o\to S_n$ be the event that $o\to x$ for some $x\in S_n$.
The following theorem, which is the site percolation version of the combination of
Theorems 6.10 and 6.14 in \cite{GG}, gives an exponential bound
on the diameter of the cluster containing the origin,  when $p$ is small.
\begin{theorem}\label{perc}
There exists a function $\phi(p)$ of $p=p(\varepsilon_0)$ such that
$$P\{o\to S_n\}\le C n^{d-1} e^{-n\phi(p)}$$
and $\lim_{p\to 0}\phi(p)=\infty$.
\end{theorem}

%GGG 04/06/2010
Let $A_x(n)$ denote the connected component of $A_x\cap \Delta_n(x)$ that contains $x$ and set 
\[q_n=P\{o\to S_n\}.\]
 The proof of Theorem
\ref{perc} will proceed by showing some (approximate) subadditivity
properties of $q_n$. We thus recall Fekete's subadditivity lemma (\cite{Fe}):

%We will
%show that $\log q_n$ satisfies certain subadditive properties and then
%prove the theorem by a lemma:
%
\begin{lemma}
	\label{lem-subadd}
	%[Theorem I\negthinspace I.2 in$\cite{GG}$]
If a sequence of finite numbers $\{b_k: k\ge 1\}$ is subadditive, that is,
$b_{m+n}\le b_m +b_n \mbox{   for all m,n}$,
then $\lim_{k\to\infty}b_k/k=\inf_{k\in\mathbb{N}} b_k/k$.
\end{lemma}
\noindent{\it Proof of Theorem \ref{perc}:}
%GGG0827
We follow the proof given by Grimmett in \cite{GG} in the bond percolation case. 
By the BK inequality (\cite{GG}, pg. 38),
$$q_{m+n}\le \sum_{x\in S_m} P\{o\to x\} P\{x\to x+S_n\}.$$
But $P\{o\to x\}\le q_m$ for $x\in S_m$
and $P\{x\to x+S_n\}=q_n$ by translation invariance. Hence we get
\begin{equation}\label{sub1}
q_{m+n}\le |S_m|q_m q_n.
\end{equation}
By exchanging $m$ and $n$ in (\ref{sub1}),
\begin{equation}\label{sub11}
q_{m+n}\le |S_{m\wedge n}|q_m q_n.
\end{equation}

On the other hand, let $U_x$ be the event that 
%GGG 04/06/2010 delete $x\in\overline{\Delta_m\cap A_o}$
$x\in \overline{A_o(m)}$
 and let $V_x$ be the event that 
%GGG 04/06/2010 delete $\overline{(x+\Delta_n)\cap A_x}\cap S_{m+n}\neq \emptyset$.
$\overline{A_x(n)}\cap S_{m+n}\neq\emptyset$. 
%GGG 03/11/2010 \cap(x+S_n) is deleted: \overline{(x+\Delta_n)\cap A_x}\cap(x+S_n)\cap ---> S_{m+n}>overline{(x+\Delta_n)\cap A_x}\cap S_{m+n}
We use the FKG inequality (\cite{GG}, pg. 34) to find that
$$q_{m+n}\ge P\{U_x\}P\{V_x\}\quad \mbox{ for any $x\in S_m$}.$$
However, $\sum_{x\in S_m}P\{U_x\}\ge q_m$, which implies that
$$\max_{x\in S_m}P\{U_x\}\ge \frac{q_m}{|S_m|} .$$
%GGG 03/06/2010 delte Let $\gamma_n$ be the probability that there exists an open path from $o$ to a face, say
%$S_n \cap \{x: x_1=n\}$, of $S_n$, 
Let $\gamma_n=P\{\overline{A_o(n)}\cap\{x:x_1=n\}\neq\emptyset\}$, then $P\{V_x\}\ge\gamma_n$.
%GGG 03/11/2010 $P\{V_x\}=\gamma_n$-->$P\{V_x\}\ge \gamma_n$
 Moreover, $\gamma_n\le q_n\le 2d\gamma_n$.
Hence
\begin{equation*}
q_{m+n}\ge \frac{q_m q_n}{2d |S_m|},
\end{equation*}
and then
\begin{equation}\label{sub2}
q_{m+n}\ge \frac{q_m q_n}{2d |S_{m\wedge n}|}.
\end{equation}

Note that $|S_m|\le C_d m^{d-1}$. Letting
$$b_k =\log q_k +\log C_d +(d-1)\log (2k),$$
one checks
using (\ref{sub11})
that the sequence $\{b_k\}$ is subadditive.
Similarly by (\ref{sub2}),  $\{-\log q_k +\log (2d C_d) +(d-1)\log (2k)\}$
is subadditive.
Thus, using Lemma \ref{lem-subadd},
$$\phi (p):=-\lim_{k\to\infty}\frac{1}{k}\log q_k$$
exists and
\begin{equation}\label{sub3}
\log q_k +\log C_d +(d-1)\log (2k)\ge -k\phi (p)
\end{equation}
\begin{equation}\label{sub4}
-\log q_k +\log (2d C_d) +(d-1)\log (2k)\ge k\phi(p).
\end{equation}
The first part of the theorem follows simply from (\ref{sub4}), and the second
by noting that with
$p\downarrow 0$ in (\ref{sub3}) we have
$q_k \downarrow 0$ and then  $\phi (p)\to \infty$. \qed
\begin{remark}
	It follows from Theorem \ref{perc}
	that
\begin{equation}\label{lo}
P\{l_o \ge n\} \le P\{o\to S_{\lfloor n/2d\rfloor}\}\le C e^{\phi(p)} n^{d-1} e^{-n\phi(p)/2d}.
\end{equation}
With (\ref{lo}) and the Borel-Cantelli lemma one concludes
that P-almost surely, $l_x\le N$ is true for all
$x\in \Delta_N$ when $N$ is sufficiently large
%OOO
and $p$ is such that $\phi(p)>0$.
Hence the inequality (\ref{f7}) holds for all $x\in\Delta_N$ when $N$
is large.
\end{remark}

\noindent{\it Proof of (\ref{Phi0}) under the assumption of Theorem \ref{CLT2}(i):}
By H\"older's inequality,
$$\frac{1}{|\Delta_N|}\sum_{y\in \partial A_x \cap \Delta_N} \Phi(y)\le
\lVert         \Phi_N 1_{\partial A_x}\rVert _{\Delta_N, \beta}
\big(\frac{|\partial A_x|}{|\Delta_N|}\big)^{1-1/\beta} ,$$
so when $N$ is large enough we have by Lemma \ref{Phicontrol} that for any $x\in A\cap \Delta_N$,
\begin{equation}\label{f8}
\Phi_N (x)\le \xi_0^{-l_x} |\partial A_x|^{1-1/\beta	} |\Delta_N |^{1/\beta}
\lVert         \Phi_N 1_{\partial A_x}\rVert _{\Delta_N, \beta}.
\end{equation}

Hence for any $\alpha \in (1, \beta)$,
\begin{align*}
&\lVert   \Phi_N 1_A\rVert _{\Delta_N,\alpha}^\alpha\\
&\le 
 \frac{1}{|\Delta_N|}\sum_{x\in A\cap\Delta_N}\big(\xi_0^{-l_x} |\partial A_x|^{1-1/\beta} 
 |\Delta_N |^{1/\beta}\lVert   \Phi_N 1_{\partial A_x}\rVert _{\Delta_N, \beta} \big)^\alpha\\
&\le 
\left[\frac{1}{|\Delta_N|}\sum_{x\in A\cap\Delta_N}
\big(\xi_0^{-l_x} |\partial A_x|^{1-1/\beta}|A_x|^{1/\beta})^{\alpha (\beta/\alpha)'}\right]
^{1-\alpha/\beta}\\
&\qquad\times \left[\frac{1}{|\Delta_N|}\sum_{x\in A\cap\Delta_N}
\big(\frac{|\Delta_N|^{1/\beta}\lVert   \Phi_N 1_{\partial A_x}\rVert _{\Delta_N, \beta}}
{|A_x|^{1/\beta}}\big)^{\beta}
\right]^{\alpha/\beta}\\
&=
 \left[\frac{1}{|\Delta_N|}\sum_{x\in A\cap\Delta_N}
\big(\xi_0^{-l_x} |\partial A_x|^{1-1/\beta}|A_x|^{1/\beta}\big)^{\alpha\beta/(\beta -\alpha)}\right]
^{1-\alpha/\beta}\\
&\qquad\times
\left(\sum_{x\in A\cap\Delta_N}
\frac{\lVert  \Phi_N 1_{\partial A_x}\rVert _{\Delta_N, \beta}^{\beta}}{|A_x|}
\right)^{\alpha/\beta},
\end{align*}

where we used (\ref{f8}) in the first inequality and H\"older's inequality in the second.

Observe that
\begin{equation}\label{Phiaverage}
\sum_{x\in A\cap\Delta_N} \dfrac{\lVert         \Phi_N 1_{\partial A_x}\rVert _{\Delta_N, \beta}^{\beta}}{|A_x|}
\le \sum_{i=1}^n \lVert         \Phi_N 1_{\partial A_i}\rVert _{\Delta_N, \beta}^{\beta}
\le 2d \lVert         \Phi_N 1_{\partial A}\rVert _{\Delta_N, \beta}^{\beta}
\le C\varepsilon_0^{-\beta},
\end{equation}
where $A_1, \cdots, A_n$ are different clusters that intersect with $\Delta_N$.
On the other hand, the multidimensional ergodic theorem gives
\begin{align}\label{percergodic}
&\lim_{N\to\infty}\frac{1}{|\Delta_N|}\sum_{x\in A\cap\Delta_N}
\big(\xi_0^{-l_x} |\partial A_x|^{1-1/\beta}|A_x|^{1/\beta}\big)^{\alpha\beta/(\beta -\alpha)}\nonumber\\
&= E \big(\xi_0^{-l_o} |\partial A_o|^{1-1/\beta}|A_o|^{1/\beta}\big)^{\alpha\beta/(\beta -\alpha)}
\le
C E  \big(\xi_0^{-l_o} l_o^d\big)^{\alpha\beta/(\beta -\alpha)} \quad \mbox{P-a.s.,}
\end{align}
which by (\ref{lo}) is finite when $\varepsilon_0$ is small.
\qed

%OOO
%GGG
\section{Transience in general ergodic environments}\label{SeTran}
In this section we will prove (ii) of Theorem
\ref{CLT1} by an argument similar to
%GGG
\cite{ZO}. The main differences in our method are that we use a stronger control of the hitting time
(Lemma \ref{tau}), and that we apply a mean value
inequality (Theorem \ref{mvi}) instead of the discrete Harnack
inequality used in \cite{ZO}.

\begin{lemma}\label{tau}
%OOO1
Let $\{X_n\}$ be a
random walk in a balanced environment
$\omega$ such that $\omega(x,o)=0$ for all $x$. For
any $r>0$, define $\tau=\tau(r)=
\inf\{n: |X_n|>r\}$. Then $E_\omega^o \tau\le (r+1)^2$.
\end{lemma}
\pf Observe that $\{|X_n|^2-n\}$ is a (quenched) martingale with respect to
$\{\mathcal{F}_n=\sigma(X_1,\cdots,X_n)\}$.
%GGG02/19/2010 The following is deleted
%Consider the $l^1$-distance from the origin. Observe that
%$P_\omega^o\{|X_{n+1}|_1-|X_n|_1=\pm 1\}=1$, and
%\begin{equation*}
%P_\omega^o\{|X_{n+1}|_1-|X_n|_1=1 | X_n\}\left\{
%                                 \begin{array}{rl}
%                                  =1/2 & \text{ if } X_n\notin \bigcup_{i=1}^d F_i,\\
%                                  >1/2 & \text{ if } X_n\in \bigcup_{i=1}^d F_i,
%                                 \end{array}
%                                \right.
%\end{equation*}
%where $F_i=\{x=(x_1,\cdots, x_n): x_i=0\}$. We have
%$$E_\omega^o(|X_{n+1}|_1^2|X_n)\ge |X_n|_1^2+1.$$
%Hence $\{|X_n|_1^2-n\}$ is a (quenched)
%submartingale with respect to $\{\mathcal{F}_n=\sigma(X_1,\cdots,X_n)\}$.
Thus by optional stopping, $0= E_\omega^o[|X_\tau|^2-\tau]\le (r+1)^2-E_\omega^o\tau$. \qed
%GGG08/25/2011: |X_\tau|_1^2 should be |X_\tau|^2

%GGG02/15/2010 added"To prove Theorem \ref{CLT1}(ii)"
To prove Theorem \ref{CLT1}(ii), we shall make use of the following mean-value inequality,
which is a modification of Theorem 3.1 in \cite{KT}. Let $B_{r}(z)=\{x\in\mathbb{Z}^d: |x-z| <r\}$. 
%be the $l^2$-ball centered at $z$. 
We shall also write $B_r (o)$  as $B_r$.
\begin{theorem}\label{mvi}
For any function $u$ on $\bar B_R (x_0)$ such that
$$L_{\omega} u =0 , \quad x \in B_R (x_0)$$
and any $\sigma\in (0,1)$, $0<p\le d$,  we have
\[
\max_{B_{\sigma R}(x_0)}u\le C \lVert
\frac{u^+}{\varepsilon^{d/p}}\rVert _{B_R (x_0), p},
\]
where $C$ depends on $\sigma$, $p$ and $d$.
\end{theorem}
We postpone the proof of Theorem \ref{mvi} to 
%GGG3/31/2010 the end of the section, 
the next section, and bring now the

\noindent{\it Proof of Theorem \ref{CLT1}(ii):}
%GGG02/15/2010
%OOO1
%GGG 3/11/2010
As mentioned in Section \ref{introduction},  the transience of the random walk would not change
if we considered the walk restricted
to its jump times. That is, the transience or recurrence of the random walk
 in an environment $\omega$ is the same as in an environment
$\tilde \omega$,
where $\tilde \omega$ is defined by
$\tilde{\omega}(x,e)=\omega (x,e)/(1-\omega(x,o))$. Therefore, in the sequel we assume $\omega(x,o)=0$
for all $x$ and almost all $\omega$.\\

Let $K$ be any constant$\ge 3$.
We denote $B_{K^i}(x)$ by $B^i(x)$ and define
$\tau_i:=\inf \{n: |X_n|> K^i\}$.
%GGG 03/11/2010 delete: We may sometimes write $X_{\tau_i} (\omega)$ instead of
%$X_{\tau_i}$ in case we want to emphasize that its randomness is given by $\omega$.
Our approach is to bound the (annealed) expected number of visits to the origin by
the walk; this requires some a-priori bounds on the moments of
$\varepsilon(o)^{-1}$.\\

For any $z\in\partial B^i$ , $y\in B^{i-1}$, noting that $E_\omega^x (\mbox{\# visits at $y$ before $\tau_{i+2}$}):=v(x)$ satisfies $L_\omega v (x)=0$ for
$x\in B^{i+2}\setminus\{y\}$, we have that for $p\in(0,d]$,
\begin{align}\label{f11}
& E_{\theta^{y}\omega}^{z} (\mbox{ \# visits at $o$ before $\tau_{i+1}$})
\nonumber\\
&\le 
 E_{\omega}^{z+y}(\mbox{\# visits at $y$ before $\tau_{i+2}$})\nonumber\\
&\le 
 \max_{x\in B^{i-1}(z)}E_{\omega}^{x}(\mbox{\# visits at $y$ before $\tau_{i+2}$})\nonumber\\
&\le 
 C\Bnorm{\frac{E_{\omega}^{x}(
               \mbox{\# visits at $y$ before $\tau_{i+2}$})}{\varepsilon_\omega (x)^{d/p}}
              } _{B_{2K^{i-1}}(z), p}\nonumber\\
&\le 
 C\Bnorm{ \frac{E_{\omega}^{x}
              (\mbox{\# visits at $y$ before $\tau_{i+2}$})}{\varepsilon_\omega (x)^{d/p}}
              } _{B^{i+2}, p},
\end{align}
where we used Theorem \ref{mvi} in the third inequality. Take $p=d/q$ (without loss of
generality, we always assume that $q< d$).
%GGG(02/08/2010) "where we used...in the second-> in the third"
Then by (\ref{f11}) and Lemma \ref{tau},
\begin{align}\label{f12}
& \sum_{y\in B^{i-1}}E_{\theta^{y}\omega}^o
(\mbox{ \# visits at $o$ in $[\tau_i,\tau_{i+1})$})\nonumber\\
&\le 
 C\sum_{y\in B^{i-1}}  \left[
\frac{1}{|B^{i+2}|}\sum_{x\in B^{i+2}}
\frac{
E_{\omega}^{x}
(\mbox{\# visits at $y$ before $\tau_{i+2}$})^{d/q}}
{\varepsilon_\omega (x)^d}
\right]^{q/d}\nonumber\\
&\le 
 C K^{-iq}\sum_{y\in B^{i-1}}\sum_{x\in B^{i+2}}
    \frac{E_{\omega}^{x}(\mbox{\# visits at $y$ before $\tau_{i+2}$})}
    {\varepsilon_\omega (x)^q}
        \nonumber\\
&= 
 C K^{-iq}\sum_{x\in B^{i+2}}
    \frac{E_{\omega}^{x}(\mbox{\# visits at $B^{i-1}$ before $\tau_{i+2}$})}
    {\varepsilon_\omega (x)^q}\nonumber\\
&\le 
 C K^{-iq}\sum_{x\in B^{i+2}}\frac{E_{\omega}^{x}\tau_{i+2}}{\varepsilon_\omega (x)^q}\nonumber\\
&\le 
 C K^{(2-q)i} \sum_{x\in B^{i+2}}       \varepsilon_\omega(x)^{-q}.
\end{align}
Taking expectations and using translation invariance we have
\begin{equation*}
%GGG delete \label{f13}
 \mathbb{E}^o(\mbox{\# visits at $o$ in $[\tau_i, \tau_{i+1})$})
\le C K^{(2-q)i} E \varepsilon^{-q}.
\end{equation*}
Therefore, if $E \varepsilon^{-q}<\infty$ for some $q>2$ , then
\begin{equation*}
\mathbb{E}^o (\mbox{\#
visits at $o$})\le C E \varepsilon^{-q}
                     \sum_{i=1}^{\infty}K^{(2-q)i}<\infty .
\end{equation*}

This proves Theorem \ref{CLT1}(ii) for $\{\Omega, P\}$ such that $\omega (x,o)=0$ for all $x$
and almost all $\omega$.
%OOO1
As mentioned earlier,
the general case follows by
replacing
$\varepsilon$ with $\varepsilon/(1-\omega(o,o))$.\qed\\

\begin{remark}
%GGG
It is natural to expect
 that arguments similar to the proof of the invariance principle also work for
proving the transience in the i.i.d. case. Namely, one may hope to control
$P_\omega^x\{\mbox{visit $o$ in $[\tau_i,\tau_{i+1})$}\}$ using some mean value inequality (like Theorem \ref{mvi}), and to use percolation arguments to handle 
``bad sites''
where the ellipticity constant $\varepsilon$ is small.

This suggests considering walks that jump from bad sites to 
good sites. In \cite{KT2}, Kuo and Trudinger
proved maximum principle and mean value inequality for balanced operators 
in general meshes, which may be applied to balanced walks with possibly
big jumps. However, their estimates 
in the presence of a small ellipticity
constant are not strong enough. To overcome this issue,
%when the ellipticity constant behaves rather badly, they don't give a good control. 
%But the transience is affected by bad ellipticity coefficients. For
%example, in $\mathbb{Z}^3$, suppose $\omega(x,e_i)\thickapprox 0$ for $i=1,2$. Although the walker may stay in the $e_3$ direction for a long time, due to ellipticity, he will walk to the $e_1$ direction, and then the $e_2$ direction in finitely many steps. So intuitively, it still behaves like a simple random walk. In the next section 
we will prove a modified maximum principle that involves only 
big exit probabilities, and then use it to prove the transience 
in the i.i.d case with no moment assumptions.
\end{remark}

\section{Transience in i.i.d. environment}\label{SeTriid}
In this section we prove a modified maximum principle for 
balanced environments.
We then prove Theorem 2(ii) using the corresponding mean 
value inequality (Theorem \ref{mvi2})
and percolation arguments.

\subsection{Balanced difference operators}\label{SeTriid1}
Following \cite{KT2}, we introduce general balanced difference operators.
Let $a$ be a nonnegative function on $\mathbb{Z}^d\times \mathbb{Z}^d$ 
such that for any $x$,
$a(x,y)> 0$ 
for only finitely many $y$. 
Define the linear operator $L_a$ acting on the set of functions
on $\mathbb{Z}^d$ by
\[L_a f(x)=\sum_y a(x,y)(f(y)-f(x)).\]
We say that $L_a$ is \textit{balanced} if
\begin{equation}\label{e7}
\sum_y a(x,y)(y-x)=0.
\end{equation}
Throughout this section we always assume that $L_a$ is a balanced probability operator, that is,
\[\sum_y a(x,y)=1.\]
For any finite subset $E\subset\mathbb{Z}^d$, define its boundary
\[E^b=E^b(a)=\{y\notin E: a(x,y)>0 \text{ for some } x\in E
      \}, 
\] 
and set 
\begin{equation}
	\label{eq-ofernew}
	\tilde{E}=E\cup E^b.
\end{equation}
Define the upper contact set of $u$ at $x\in E$ as
$$I_u(x)=I_u(x,E,a)=\{s\in\mathbb{R}^d: u(x)-s\cdot x\ge u(z)-s\cdot z \text{ for all }z\in\tilde{E}\}.$$
Set 
$$h_x=h_x(a)=\max_{y: a(x,y)>0}\abs{x-y}.$$

The following lemma is useful in the proofs of various mean 
value inequalities. It is similar to Theorem 2.2 in \cite{KT2}, 
except that the proof in \cite{KT2} contains several unclear
 passages, e.g. in the inequality above (2.23) in \cite{KT2},
and so we provide a complete proof.
Throughout, we set $u^+=u\vee 0$. 
\begin{lemma}\label{mvilemma}
%GGG0827
Fix $R>0$. Let $\eta(x)=\eta_R(x):=(1-\abs{x}^2/R^2)^\beta 1_{|x|<R}$ be a
function on $\mathbb{R}^d$.
For any function $u$ on $B_R$ such that 
$L_a u=0$ in $B_R$ and any $\beta\ge 2$, 
we let $v=\eta u^+$.
Then for any $x\in B_R$ with $I_v (x)=I_v(x, B_R, a)\neq\emptyset$,
\[L_a v(x)\ge -C(\beta) \eta^{1-2/\beta}R^{-2} h_x^2 u^+,\]
where $C(\beta)$ is a constant that depends only on $\beta$.
\end{lemma}
\pf 
%GGG
We only need to consider the nontrivial case that $v\not\equiv 0$.
For $s=s(x)\in I_v(x)\neq \emptyset$,  recalling the definition of $I_v$ one has that
$$|s|\le 2v(x)/(R-|x|).$$ 
%GGG 
Note that $I_v(x)\neq\emptyset$ implies $u(x)> 0$.
If further $R^2-|x|^2\ge 4R \abs{x-y}$ , 
computations as in \cite[pg. 426]{KT2}
reveal that
% (\cite{KT2}, pg.426):
\begin{align}
 2^{-\beta}
 &\le
 \frac{\eta(y)}{\eta(x)}\le 2^\beta,\label{first}\\
 \abs{\eta(x)-\eta(y)}
 &\le
  \beta 2^\beta R^{-1} \eta(x)^{1-1/\beta}|x-y|,\label{second}\\
 \abs{\eta(x)-\eta(y)-\nabla\eta(x)(x-y)}
 &\le
 \beta (\beta-1)2^{\beta} R^{-2}\eta(x)^{1-2/\beta}|x-y|^2,\label{third}\\
%GGG082311
 |s|
 &\le 
 4 \eta^{1-1/\beta}R^{-1}u,\label{fourth}
 \end{align}
%GGG0827
where $\nabla\eta$ is the gradient of $\eta$.
 Following \cite{KT2}, we set $w(z)=v(z)-s\cdot (z-x)$. 
By the definition of $s$, we have
 $w(x)\ge w(z)$ for all $z\in \tilde E$ and
\begin{align}\label{e4}
&\sum_y a(x,y)\big(v(x)-v(y)\big)\nonumber\\
&\stackrel{(\ref{e7})}{=}
\sum_y a(x,y)\big(w(x)-w(y)\big)\nonumber\\
&\stackrel{(\ref{first}),\, w(x)\geq w(y)}{\le} 
 2^\beta \sum_y a(x,y) \frac{\eta (x)}{\eta (y)}\big(w(x)-w(y)\big)\nonumber\\
&= 
 2^\beta \sum_y a(x,y) \Big[\frac{\eta(x)}{\eta(y)}\big(v(x)-v(y)\big)+\frac{\eta(x)}{\eta(y)}s (y-x)\Big].
\end{align}
Consider first $x$ such that $R^2-|x|^2\ge 4Rh_x$.
%GGG
Then (recalling that $u(x)> 0$ because $I_v(x)\neq \emptyset$),
\begin{align}
& \sum_y a(x,y) \frac{\eta(x)}{\eta(y)}\big(v(x)-v(y)\big)\nonumber\\
&= 
\sum_y a(x,y)\left [\eta(x)\big(u(x)-u^+(y)\big)+\big(\eta(x)-\eta(y)\big)u(x)+\frac{(\eta(x)-\eta(y))^2}{\eta(y)}u(x)\right]\nonumber\\
&\stackrel{a\geq 0}{\le}  
\eta(x)L_a u(x)+\sum_y a(x,y) \left[\big(\eta(x)-\eta(y)\big)u(x)+\frac{(\eta(x)-\eta(y))^2}{\eta(y)}u(x)\right]\nonumber\\
&\stackrel{L_au=0,\, \eqref{e7}}{=}
 \sum_y a(x,y) \left[\big(\eta(x)-\eta(y)-\nabla\eta(x)(x-y)\big)u(x)+\frac{(\eta(x)-\eta(y))^2}{\eta(y)}u(x)\right]\nonumber\\
&\le  
\beta^2 2^{3\beta +1} \eta^{1-2/\beta}R^{-2}  h_x^2 u,\label{e8}
\end{align}
where we used 
%(\ref{e7}) and $L_a u=0$ in the second equality and 
(\ref{first}), (\ref{second}), (\ref{third})
in the last inequality.
Moreover, by (\ref{e7}), (\ref{first}), (\ref{second}) and (\ref{fourth}),
\begin{align}\label{e9}
\sum_y a(x,y)\frac{\eta(x)}{\eta(y)}s (y-x)
&= 
\sum_y a(x,y) \frac{\eta(x)-\eta(y)}{\eta(y)} s\cdot (y-x)\nonumber\\
&\le  \beta 2^{2\beta+2} \eta^{1-2/\beta}R^{-2} h_x^2 u.
\end{align}
Hence, combining (\ref{e4}), (\ref{e8}) and (\ref{e9}), we conclude that 
\[
-L_a v\le \beta^2 2^{4\beta+2}\eta^{1-2/\beta} R^{-2} h_x^2 u
\]
holds in $\{x: R^2-|x|^2\ge 4Rh_x, I_v(x)\neq\emptyset\}$.

On the other hand, if $R^2-|x|^2<4Rh_x$,  then $\eta^{1/\beta}\le 4h_x/R$. Thus by the fact that
$v\ge 0$, we have $-L_a v\le 2v(x)\le 32\eta^{1-2/\beta} R^{-2} h_x^2 u$. \qed\\

\noindent{\it Proof of Theorem \ref{mvi}:}
Since $L_\omega$ is a balanced operator and $h_x=1$ in this case,
by the above lemma, 
\[L_\omega v\ge -C(\beta) \eta^{1-2/\beta}R^{-2} u\] for $x\in B_R$ such that $I_u(x)\neq\emptyset$.
Applying Theorem \ref{MP} to $v$ and taking $\beta=2d/p\ge 2$, 
we obtain
%GGG0827
\begin{align*}
\max_{B_R} v 
&\le 
C \Bnorm{ \eta^{1-2/\beta} \frac{u^+}{\varepsilon}} _{B_R, d}
=C\Bnorm{    v^{1-p/d}\frac{(u^+)^{p/d}}{\varepsilon}}_{B_R, d}\\
&\le 
C(\max_{B_R} v)^{1-p/d}\Bnorm{\frac{u^+}{\varepsilon^{d/p}}} _{B_R, p}^{p/d}.
\end{align*}
Hence 
\[
\max_{B_R} v \le C \Bnorm{\frac{u^+}{\varepsilon^{d/p}}} _{B_R, p},
\]
and then
\[
\max_{B_{\sigma R}}u\le (1-\sigma^2)^{-2d/p}\max_{B_{\sigma R}}v
\le C(\sigma, p, d) \Bnorm{\frac{u^+}{\varepsilon^{d/p}}} _{B_R, p}. \text{\qed}
\]

%The theorem follows. 

\subsection{A new maximum principle and proof of Theorem \ref{CLT2}(ii)}

For any fixed environment $\omega\in\Omega$, let $\varepsilon_0>0$ be a constant to be determined,
and define site percolation as in Section \ref{SePercE}. Recall that for $x\in\mathbb{Z}^d$, $A_x$
is the percolation cluster that contains $x$ and $l_x$ is its $l^1$-diameter.
As mentioned in the introduction, the transience would not change if we 
considered the walk restricted to its
jump times. Without loss of generality, we assume that $\omega(x,o)=0$ for all $x$, $P$-almost surely.

Recall the definition of $\square$ and $\kappa$-path for $\kappa\in\square$ in Section \ref{SePercE}. Note
that under our assumption, $\max_i \omega(x, e_i)\ge 1/2d$, so we take $\xi_0=1/2d$ in the definition of 
$\kappa$-paths.

For each $\kappa\in\square$, we pick a site $y_\kappa=y(x, \kappa)\in \partial A_x$ such that
\[d(x, y_\kappa)=\max_{\substack{y:
                  \exists \text{ $\kappa$-path in $\bar{A}_x$ }\\\text{ from $x$ to $y$}
                           }
                 } d(x,y)\] 
and let $\Lambda_x\subset \bar{A}_x$ be the union of (the points of the) $\kappa$-paths from $x$ to $y_\kappa$ over all $\kappa\in\square$.
From the definition of $y_\kappa$ one can conclude that
\begin{itemize}
\item For any $q\in\mathbb{R}^d$, we pick a $\kappa=\kappa_q\in \square$ such that
\[q_j \kappa_j\le 0 \text{ for all }j=1,\cdots, d.\]
Then $(y_\kappa-x)_j q_j\le 0$ for all $j=1,\cdots, d$. Moreover, for $i\in\{1,\cdots, d\}$, $q_i>0$ implies $y_\kappa-e_i\notin \Lambda_x$, and $q_i<0$
implies $y_\kappa+e_i\notin \Lambda_x$.
\end{itemize}

In the sequel we let $\tau_{\Lambda_x}=\inf\{n>0: X_n\notin \Lambda_x\}$ and 
\[a(x,y)=P_\omega^x \{
                      X_{\tau_{\Lambda_x}}=y
                     \}.
\]
By the fact that $X_n$ is a (quenched) 
martingale, it follows that $L_a$ is a balanced operator.

%OOO0828
For the statement of the next theorem,
recall the definition of $\tilde{E}$, see \eqref{eq-ofernew}.
\begin{theorem}\label{mp2}
Let $E\subset\mathbb{Z}^d$ be bounded. 
%Recall the definition of $\tilde{E}$ in Section \ref{SeTriid1}. 
Let $u$ be a function on $\tilde{E}$. If $L_a u(x)\ge -g(x)$ for all $x\in E$ such that $I_u(x)=I_u(x,E,a)\neq \emptyset$ ,
then
$$\max_E u\le 
  \frac{d\diam \tilde{E}}{\varepsilon_0}
  \bigg(\sum_{\substack{
                           x\in E\\I_u (x)\ne \emptyset
                       }
             }
                      \abs{g(x)(2d)^{l_x}}^d
  \bigg)^{\frac{1}{d}}
 +\max_{E^b}u  .$$
\end{theorem}
\pf Without loss of generality, assume $g\ge 0$ and
$$\max_E u=u(x_0)>\max_{E^b}u$$
for some $x_0\in E$. Otherwise, there is nothing to prove.

For $ s\in \mathbb{R}^{d} $ such that $ |s|_{\infty} \le [u(x_0)-\max_{E^b}u]/(d\diam \tilde{E}) $,
we have \[ u(x_0)-u(x) \ge s \cdot (x_0-x) \] for all $ x \in E^b $,
%GGG
 which implies that
$\max_{z\in\tilde{E}}u(z)-s\cdot z$
is achieved in $E$. Hence
  $ s \in \bigcup_{x \in E} I_{u}(x) $ and

\begin{equation}\label{e1}
 \left[-\dfrac{u(x_0)-\max_{E^b }u}{d\diam \tilde{E}}, \dfrac{u(x_0)-\max_{E^b}u}{d\diam \tilde{E}}\right]^{d} \subset \bigcup_{x\in E} I_{u}(x)  .
 \end{equation}

Further, if $s\in I_u(x)$, we set
$$w(z)=u(z)-s(z-x).$$
Then $w(z)\le w(x)$ for all $z\in \tilde{E}$ and 
\begin{equation}\label{e2}
I_u(x)=I_w(x)+s.
\end{equation}
Since for any $q\in I_w(x)$,  there is $\kappa=\kappa_q\in\square$ such that
$$q_j(x-y_\kappa)_j\ge 0 \text{ for } j=1,\cdots, d,$$ we have
$$ w(x)-w(y_\kappa\pm e_i)\ge q(x-y_\kappa\mp e_i)\ge \mp q_i.$$
Moreover, for any $i\in\{1,\cdots, d\}$, if $q_i>0$, then $y_\kappa-e_i\notin \Lambda_x$ and we have
$w(x)-w(y_\kappa-e_i)\ge |q_i|$. Similarly, if $q_i<0$, then $y_\kappa+e_i\notin \Lambda_x$ and 
$w(x)-w(y_\kappa+e_i)\ge |q_i|$. We conclude that
$$|q_i|\le \frac{\sum_{y} a(x, y)(w(x)-w(y))}{\min\limits_{\pm}\{a(x, y_\kappa\pm e_i)\}}.$$
On the other hand, from the construction of $\Lambda_x$ we obtain (note that $y_\kappa\in\partial A_x$)
$$a(x, y_\kappa\pm e_i)\ge (\frac{1}{2d})^{l_x} \varepsilon_0.$$
Hence, since $L_a$ is balanced,
$$|q_i|\le 
\frac{(2d)^{l_x}}{\varepsilon_0}\sum_y a(x,y)(w(x)-w(y))
=
\frac{(2d)^{l_x}}{\varepsilon_0}(-L_a u)
%\sum_y a(x,y)(w(x)-w(y))
\le \frac{(2d)^{l_x}}{\varepsilon_0} g$$
for all $i$. Therefore
\begin{equation}\label{e3}
I_w(x)\subset [-(2d)^{l_x}\varepsilon_0^{-1} g, (2d)^{l_x}\varepsilon_0^{-1} g]^d.            
\end{equation}
Combining (\ref{e1}), (\ref{e2}) and (\ref{e3}) we conclude that
\begin{equation*}
\left(\dfrac{u(x_0)-\max_{E^b}u}{d\diam \tilde{E}}\right)^d\le 
\sum_{\substack{
                           x\in E\\I_u (x)\ne \emptyset
                       }
             }
                      \abs{g(x)(2d)^{l_x}\varepsilon_0^{-1}}^d. 
\quad\quad\quad\quad\quad
\quad\quad\quad\quad\quad
 \mbox{\qed}
\end{equation*}

As with Theorem \ref{mvi}, we have a corresponding mean value inequality.
% for Theorem \ref{mp2}:
\begin{theorem}\label{mvi2}
For any function $u$ on $B_R$ such that
$$L_a u=0, \quad x\in B_R$$ 
and any $\sigma\in (0,1)$, $0<p\le d$, we have
\[
\max_{B_{\sigma R}} u\le C\big(\frac{\diam \tilde{B}_R}{\varepsilon_0 R}\big)^{d/p}
                           \norm{[l_x^2 (2d)^{l_x}]^{d/p}u^+}_{B_R, p},
\]
where $C$ depends on $\sigma, p$ and $d$.
\end{theorem}
\pf By the same argument as in the proof of Theorem \ref{mvi}, Lemma \ref{mvilemma} and 
Theorem \ref{mp2} implies Theorem \ref{mvi2}. \qed\\

Having established
Theorem \ref{mvi2}, we can now 
prove the transience of the random walks in balanced
i.i.d. environment with $d\geq 3$.

\noindent{\it Proof of Theorem \ref{CLT2}(ii)}:
Let $K$ be any constant $\ge 4$ and define $B^i, \tau_i$ as in Section 5.
Let $\Omega_i=\{\omega\in\Omega: l_x\le K^{i-1} \mbox{ for all $x\in B^{i+2}$}\}$. 
For any $\omega\in\Omega_i$, $z\in\partial B^i$, $y\in B^{i-1}$, noting that
$P_\omega^x \{\mbox{visit $y$ before $\tau_{i+2}$}\}:=u(x)$ satisfies
\[L_a u(x)=0\] 
%GGG
for $x\in B_{2K^{i-1}}(z)$, by similar argument as in (\ref{f11}) we have
\begin{align*}
& P_{\theta^y\omega}^z \{\mbox{visit $o$ before $\tau_{i+1}$}\}1_{\omega\in\Omega_i}\\
&\le 
\max_{x\in B^{i-1}(z)}P_\omega^x \{\mbox{visit $y$ before $\tau_{i+2}$}\}1_{\omega\in\Omega_i}\\
&\le 
C \varepsilon_0^{-d}\norm{
                                [l_x^2 (2d)^{l_x}]^d P_\omega^x \{\mbox{visit $y$ before $\tau_{i+2}$}\}                              
                               }_{B_{2K^{i-1}}(z), 1}\\
&\le 
C \varepsilon_0^{-d}\abs{B^{i+2}}^{-1}
         \sum_{x\in B^{i+2}}l_x^{2d} (2d)^{dl_x} P_\omega^x \{\mbox{visit $y$ before $\tau_{i+2}$}\} ,                            
\end{align*}
where in the second inequality, we applied Theorem \ref{mvi2} with
$p=1$ and used the fact that $\diam \tilde{B}_{2K^{i-1}}\le 3K^{i-1}$ when $\omega\in\Omega_i$. Hence
\begin{align}\label{e5}
&\sum_{y\in B^{i-1}}P_{\theta^y\omega}^o 
                     \{\mbox{visit $o$ in $[\tau_i, \tau_{i+1})$}\}
                     1_{\omega\in\Omega_i}\nonumber\\
&\le 
C \varepsilon_0^{-d}\abs{B^{i+2}}^{-1}
         \sum_{x\in B^{i+2}}l_x^{2d} (2d)^{dl_x} E_\omega^x 
         (\mbox{\# visits at $B^{i-1}$ before $\tau_{i+2}$})\nonumber\\
&\stackrel{\text{Lemma }\ref{tau}}{\le} 
C \varepsilon_0^{-d}K^{(2-d)i}\sum_{x\in B^{i+2}}l_x^{2d} (2d)^{dl_x}.                              
\end{align}

Since 
\begin{align}\label{e6}
&\sum_{y\in B^{i-1}}P_{\theta^y\omega}^o 
                     \{\mbox{visit $o$ in $[\tau_i, \tau_{i+1})$}\}\nonumber\\
&\le \sum_{y\in B^{i-1}}P_{\theta^y\omega}^o 
                     \{\mbox{visit $o$ in $[\tau_i, \tau_{i+1})$}\}1_{\omega\in\Omega_i}
                     +\abs{B^{i-1}}1_{\omega\notin\Omega_i},                     
\end{align}
taking $P$-expectations on both sides of (\ref{e6}) and using (\ref{e5}) we get
\[
\mathbb{P}^o \{\mbox{visit $o$ in $[\tau_i,\tau_{i+1})$}\}
\le
C \varepsilon_0^{-d}K^{(2-d)i}El_o^{2d}(2d)^{dl_o}+P\{\omega\notin\Omega_i\}.
\]
By (\ref{lo}), we can take $\varepsilon_0$ to be small enough such that
$El_o^{2d}(2d)^{dl_o}<\infty$ and $\sum_{i=1}^\infty P\{\omega\notin\Omega_i\}<\infty$.
Therefore when $d\ge 3$,
\[\sum_{i=1}^\infty \mathbb{P}^o \{\mbox{visit $o$ in 
$[\tau_i,\tau_{i+1})$}\}<\infty. \quad\quad\quad\quad\quad\quad
\mbox{\qed}\]

\section{Concluding remark}
%OOO0827
 While Bouchaud's trap model (see \cite{Bou,BAC}) provides an example of an (i.i.d.)
environment where local traps can destroy the invariance
principle, it is interesting to note that
a counter-example to Theorem
		\ref{CLT2} in the ergodic setup also can be written.
Namely, let $d\ge 2$, write for $x\in\mathbb{Z}^d$, $z(x)=(x_2,\cdots, x_d)\in\mathbb{Z}^{d-1}$.
Let $\{\varepsilon_z\}_{z\in\mathbb{Z}^{d-1}}$ be i.i.d random variables with support in $(0, 1/2)$ and set
\begin{equation}
\omega(x, e)=\left\{
\begin{array}{rl}
\varepsilon_{z(x)}, & \text{if } e=\pm e_1\\
(1-2\varepsilon_{z(x)})/2(d-1), & \text{else }.
\end{array}
\right.
\end{equation}
It is easy to verify that $\{X_t^n\}_{t\ge 0}$ satisfies the quenched invariance principle, but that
the limiting covariance may degenerate if the tail of $\varepsilon_z$ is heavy.


\begin{thebibliography}{99}
\bibitem{Bar04} M. Barlow, 
\emph{Random walks on supercritical percolation clusters},
 Annals Probab.  32  (2004),   3024--3084.
\bibitem{BarDe10} M. Barlow, J.-D. Deuschel,
\emph{Invariance principle for 
the random conductance model with unbounded conductances}, 
Annals Probab. 38 (2010), 234--276.
	\bibitem{BAC} G. Ben Arous, J. Cerny, \emph{Scaling limits for trap models on $\mathbb{Z}^D$}, Annals Probab. 35 (2007),  2356--2384.
\bibitem{BB} N. Berger, M. Biskup,
\emph{Quenched invariance principle for simple random walks on percolation clusters},  Probab. Theory Related Fields  137  (2007),   83--120.
%\bibitem{PB} P. Billingsley,
%\emph{Convergence of probability measures, 2nd Edition}, John Wiley, 1999.
%OOO \bibitem{BZZ} M. Bramson, O. Zeitouni,  M. P. W. Zerner,
%\emph{Shortest Spanning Trees and a Counterexample for Random Walks in
%Random Environments}, Annals Probab. 34 (2006),  821--856.
%\bibitem{Du} R. Durrett,
%\emph{Probability: Theory and Examples (Third edition)}, Brooks-Cole-Thomson Learning, 2005.
%\bibitem{DE} A. Dvoretzky, P. Erd\"os,
%\emph{Some problems on random walk in space}, Proceedings of Second Berkeley Symposium, (1951), 353--367.
%GGG0827
\bibitem{Bou}J.P.Bouchaud,
\emph{Weak ergodicity breaking and aging in disordered systems},
J. Phys. I (France) 2, 1705 (1992).
\bibitem{DFGW89} A. De Masi, P. A. Ferrari, S. Goldstein, W. D. Wick,
\emph{An invariance principle for reversible Markov processes. Applications
 to random motions in random environments}, J. Statist. Phys. 55 (1989),
787--855.
\bibitem{Fe} M. Fekete,
\emph{\"{U}ber die Verteilung der Wurzeln bei gewissen algebraischen Gleichungen mit ganzzahligen Koeffizienten}, Mathematische Zeitschrift 17 (1923), 228--249.
\bibitem{Ge} H. O. Georgii,
\emph{Gibbs Measures and Phase Transitions}, Walter de Gruyter, Berlin, 1988.
\bibitem{GG} G. Grimmett,
\emph{Percolation (Second edition)}, Springer, 1999.
\bibitem{Kr} U. Krengel,
\emph{Ergodic theorems}, Walter de Gruyter, Berlin, 1985.
\bibitem{KT} H. J. Kuo, N. S. Trudinger,
\emph{Linear elliptic difference inequalities with random coefficients}, Math. Comp. 55 (1990), 37--53.
\bibitem{KT2} H. J. Kuo, N. S. Trudinger,
\emph{Positive difference operators on general meshes},  
Duke Math. J. 83 (1996),  415--433.
\bibitem{La} G. F. Lawler,
\emph{Weak convergence of a random walk in a random environment}, 
Comm. Math. Phys. 87 (1982),  81--87.
\bibitem{Ma08} P. Mathieu, 
\emph{Quenched invariance principles for 
random walks with random
 conductances},
 J. Stat. Phys.  130  (2008),   1025--1046.
\bibitem{MR05} P. Mathieu, E. Remy, 
\emph{Isoperimetry and heat kernel decay on percolation clusters},
Annals Probab. 32 (2004),  100--128.
\bibitem{MaP07} P. Mathieu, A. Pianitski, 
\emph{Quenched invariance principles for random walks on percolation
 clusters}, Proc. R. Soc. Lond. Ser. A Math. Phys. Eng. Sci.  463  (2007),
   2287--2307.
\bibitem{PV} G. Papanicolaou, S.R.S.Varadhan,
\emph{Diffusions with random coefficients}, Statistics and probability: essays in honor of C. R. Rao,  North-Holland, Amsterdam, 1982, pp. 547--552
\bibitem{Sz} A. S. Sznitman,
\emph{Lectures on random motions in random media}, In DMV seminar 32, Birkhauser, Basel, 2002.
\bibitem{SS05} V. Sidoravicius, A. S. Sznitman,
\emph{Quenched invariance principles for walks on clusters of percolation or
 among random conductances},
 Probab. Theory Related Fields  129  (2004),   219--244.
\bibitem{ZO} O. Zeitouni,
\emph{Random walks in random environment}, XXXI Summer school in Probability, St. Flour (2001). Lecture notes in Math. 1837, Berlin:Springer, 2004,  193--312
\bibitem{Zrev} O. Zeitouni,
\emph{Random walks in random environments},
J. Phys. A  39  (2006), R433--R464.
\end{thebibliography}
\end{document}